\tolerance=10000
\magnification=1200
\raggedbottom

\baselineskip=15pt
\parskip=1\jot

\def\sk{\vskip 3\jot}

\def\heading#1{\vskip3\jot{\noindent\bf #1}}
\def\label#1{{\noindent\it #1}}


\def\ref#1;#2;#3;#4;#5.{\item{[#1]} #2,#3,{\it #4},#5.}
\def\refinbook#1;#2;#3;#4;#5;#6.{\item{[#1]} #2, #3, #4, {\it #5},#6.} 
\def\refbook#1;#2;#3;#4.{\item{[#1]} #2,{\it #3},#4.}


\def\({\bigl(}
\def\){\bigr)}

\def\ze{\zeta}

{
\pageno=0
\nopagenumbers
\rightline{\tt circle.limit.tex}
\vskip1in

\centerline{\bf Counting the Angels and Devils in Escher's {\it Circle Limit IV}}
\vskip0.5in

\centerline{John Choi}
\centerline{\tt aeatda@gmail.com}
\sk
\centerline{Yonghyeon-ro 10, Muwon Village 506-1504}
\centerline{Goyang, Gyeonggi Province}
\centerline{Republic of Korea 412-724}
\sk\sk

\centerline{Nicholas Pippenger}
\centerline{\tt njp@math.hmc.edu}
\sk

\centerline{Department of Mathematics}
\centerline{Harvey Mudd College}
\centerline{301 Platt Boulevard}
\centerline{Claremont, CA 91711}
\vskip0.5in

\noindent{\bf Abstract:}
We derive the rational generating function that enumerates the angels and devils in M.~C. Escher's {\it Circle Limit IV} according to their combinatorial distance from the six creatures whose feet meet at the center of the disk.
This result shows that the base of the exponential rate of growth is $1.582\ldots$
(the largest root of the polynomial $1 - z^2 - 2z^3 - z^4 + z^6$).
\vfill\eject
}

\heading{1. Introduction}

M.~C. Escher's {\it Circle Limit IV\/} (1960) is one of the artist's most famous works.
Its mathematical aspects have already been studied by Coxeter [C2], Bennett [B] and Dunham [D].
In it, white angels and black devils alternate to form a tessellation of the hyperbolic plane by triangular faces. 
Each creature (angel or devil) is represented by a face, with one vertex at the feet and two vertices at the wing tips, and with one edge at the head and two edges at the legs.
Each edge separates an angel and a devil.
At each vertex of the tessellation, six pairs of feet (three of angels and three of devils) or eight wing tips (four of angels and four of devils) meet. 
The tessellation is displayed with six pairs of feet meeting at the center of a Poincar\'{e} disk model of the hyperbolic plane.
The tessellation is isohedral (the symmetry group acts transitively on the faces, if one ignores the distinction between angles and devils), but not isogonal (some vertices have degree six, while others have degree eight) and not isotoxal (some edges join two vertices of degree eight, while others join a vertex of degree eight with a vertex of degree six).

Choi [C1] refers to this tessellation as $\{3, [6,8,8]\}$ (a tessellation by faces of degree three, in which each triangle meets one vertex of degree six, and two vertices of degree eight).
The vertices may be partitioned into generations, with the $n$-th generation comprising the vertices
at distance $n$ from the single vertex of degree six at the center
(where the distance between two vertices is the length of a shortest path between them, and the length of a path is the number of edges it contains).
Choi [C1] gives the generating function that enumerates vertices of this tessellation according to their generation (that is, in which the coefficient of $z^n$ is the number of vertices in generation $n$). 
This generating function is $(1+4z+10z^2+4z^3+z^4)/(1-2z-2z^2-2z^3+z^4) = 1+6z+24z^2+66z^3+192z^4+\cdots\,$,
where we have included enough terms to show that the sequence of coefficients of the series in parentheses is not in the 
{\it OEIS\/} (the {\it On-Line Encyclopedia of Integer Sequences}, at {\tt http://oeis.org})
as of this writing.
Since this generating function is rational with simple poles, the coefficients grow exponentially at a rate given by the reciprocal of the pole closest to the origin (which, since the denominator 
polynomial $1-2z-2z^2-2z^3+z^4$ is self-reciprocal, is the pole furthest from the origin).
This rate (that is, the base of the exponential)  is $1/2 + \sqrt{5}/2 + \sqrt{1/2 + \sqrt{5}/2} = 2.890\ldots\,$.
It is natural to expect that this growth rate is intermediate between that of the regular tessellation
$\{3,6\}$ (with six triangles meeting at every vertex) and that of $\{3,8\}$ (with eight triangle meeting at every vertex).
The former is a tessellation of the Euclidean plane, and therefore has linear growth
(exponential rate $1 < 2.890\ldots\,$, as expected).
Paul and Pippenger [P] have determined the generating functions for all regular tessellation of the hyperbolic plane; for $\{3,8\}$ the generating function is $(1+4z+z^2)/(1-4z+z^2)$ and the exponential rate is $2+\sqrt{3} = 3.732\ldots > 2.890\ldots\,$, as expected.

Our goal in this paper is to conduct an analogous enumeration of the faces of the tessellation 
$\{3, [6,8,8]\}$ (so that we are counting the angels and devils themselves, rather than the meeting places of their feet and wing tips.
To do this we consider the dual tessellation $\{[6,8,8], 3\}$, constructed by taking a vertex corresponding to each face in $\{3, [6,8,8]\}$ and a face corresponding to each vertex in $\{3, [6,8,8]\}$,
with an edge joining each pair of vertices in $\{[6,8,8], 3\}$ that correspond to adjacent faces in 
$\{3, [6,8,8]\}$.
In this dual tessellation, every vertex has degree three, and every vertex is the meeting place of one hexagon (face of degree six) and two octagons (faces of degree eight).
This tessellation is isogonal, but not isohedral and not isotoxal.
We are interested in enumerating the vertices of this dual tessellation.
We again partition the vertices into generation, now according to the length of their shortest path to any of the six vertices of the central hexagon (corresponding to the six creatures whose feet meet at the center of the original tessellation).
We shall show that the generating function for this dual tessellation is
$6(1+z+z^2+z^3) / (1-z^2-2z^3-z^4+z^6) = 6(1+z+2z^2+4z^3+5z^4+9z^5+14z^6+22z^7+\cdots)$.
(The sequence of coefficients of the series in parentheses is not in the{\it OEIS\/}).)
The exponential growth rate is $1.582\ldots\,$.
It is natural to expect that this growth rate is intermediate between that of the regular tessellation
$\{6,3\}$ (with three hexagons meeting at every vertex) and that of $\{8,3\}$ (with three octagons meeting at every vertex).
The former is a tessellation of the Euclidean plane, and therefore has linear growth
(exponential rate $1 < 1.582\ldots\,$, as expected).
Paul and Pippenger [P] have shown that for $\{8,3\}$ the generating function is 
$(1+z)(1+z+z^2+z^3+z^4)/(1-6z-6z^2-6z^3+z^4)$ and the exponential rate is 
$1/4 + \sqrt{13}/4 + \sqrt{\sqrt{13}/8 - 1/8} = 1.722\ldots > 1.582\ldots\,$, as expected.

The polynomial $1 - z^2 - 2z^3 - z^4 + z^6$ that determines the growth rate of vertices in $\{[6,8,8], 3\}$ is solvable (because it is self-reciprocal).
The growth rate is the larger root of $z + 1/z = \ze$, where 
$$\ze = {(9 + i\sqrt{111})^{1/3} \over 3^{2/3}} + {4\over 3^{1/3} \; (9 + i\sqrt{111})^{1/3}}$$
is one of the roots of the polynomial $\ze^3 - 4\ze - 2$.
(The other two roots of $\ze^3 - 4\ze - 2$ give rise to complex roots of $1 - z^2 - 2z^3 - z^4 + z^6$.)

The theme of alternating angels and devils was clearly a favorite of Escher's; almost twenty years earlier he had constructed analogous tessellations of the Euclidean plane (1941)
and the sphere (1942) (see Coxeter [C2] and Dunham [D]).
The analogous enumerations for these tessellations are much easier than those  done 
for {\it Circle Limit IV} by Choi [C1] and here, and we merely summarize them here.
The tessellation of the Euclidean plane may be described as $\{3, [4, 8, 8]\}$ (four pairs of feet or eight wing tips meet at each vertex); the generating function enumerating its vertices (with the origin at a meeting of feet, as for {\it Circle Limit IV\/}) is 
$(1 + 2 z + 9 z^2 - 4 z^3)/(1-z)^2 = 1 + 4z + 16z^2 + 24z^3 + 32z^4 + \cdots\,$; the growth rate is of course linear.
(This sequence is not in the {\it OEIS}.)
The creatures,  corresponding to the vertices of the dual tessellation $\{[4,8,8], 3\}$ are enumerated by the generating function $4(1+z^2)/(1+z+z^2)(1-z)^2 = 
4(1 + z + 2z^2 + 3z^3 + 3z^4 + 4z^5 + 5z^6 + 5z^7 + \cdots\,$; all the roots of the denominator are of absolute value unity, so the growth rate is again linear.
(The sequence of coefficients of the series in parentheses is A004396 in the the {\it OEIS}.)
The tessellation of the sphere is $\{3, [4, 6, 6]\}$ (four pairs of feet or six wing tips meet at each vertex), with vertices enumerated by the polynomial $1+4z+8z^+z^3$; the vertices of the dual
$\{[4,6,6]\}$ are enumerated by $4(1+z+2z^2+z^3+z^4)$.
\sk

\heading{2. Enumeration}

From this point onward, all references to vertices, edges and faces refer to the dual tessellation 
$\{[6,8,8], 3\}$.
We begin by classify the edges and defining some relations among vertices.
An edge that joins a vertex $v$ in generation $n$ and a vertex $w$ in generation $n+1$ will be called a {\it filial\/} edge; $v$ will be called a {\it parent\/} of $w$, and $w$ will be called a {\it child\/} of $v$.
An edge that joins two vertices in the same generation will be called a {\it fraternal\/} edge, and the joined vertices will be called {\it siblings\/} of each other.

Each of the six vertices in generation zero has two siblings and one child.
We shall call these $G$-vertices.
The vertices in subsequent generations  are of three types.
Those that have one parent and two children will be called  type-I vertices;
those with two parents and one child will be called type-II vertices; and
those with one parent, one sibling and one child will be called type-III vertices.

A hexagon (or octagon) that has two type-III sibling vertices  in generation $n$ and two type-III sibling vertices in generation 
$n+2$ (or $n+3$)
will be called a {\it flat\/} hexagon (or octagon).
The vertices in generation $n$ will be called the {\it floor\/} of the face, and those in generation
$n+2$ (or $n+3$) will be called the {\it roof}.
A hexagon (or octagon) that has a type-I vertex in generation $n$ and a type-II vertex in generation 
$n+3$ (or $n+4$)
will be called a {\it sharp\/} hexagon (or octagon).
The vertex in generation $n$ will be called the {\it floor\/} of the face, and that in generation
$n+3$ (or $n+4$) will be called the {\it roof}.

Type-III sibling vertices that form the roof of a hexagon and the floor of an octagon will
be called $E$-vertices; those that form the roof of an octagon and the floor of a hexagon will be called $F$-vertices.
Type-II vertices that form the roof of a hexagon will be called $C$-vertices;
those that form the roof of an octagon will be called $D$-vertices.
Type-I vertices that form the floor of a hexagon will be called $A$-vertices; those that form the 
floor of an octagon will be called $B$-vertices.

Let $G(z)$ denote the generating function for $G$-vertices.
We have
$$G(z) = 6. \eqno(1)$$

Let $E(z)$ and $F(z)$ denote the generating functions for $E$-vertices and $F$-vertices, respectively.
Each of the six pairs of consecutive $G$-vertices in generation $0$ gives rise to a pair of $F$-vertices across the flat octagon in generation $3$.
Each pair of $F$-vertices in generation $n$ gives rise to a pair of $E$-vertices across the  flat hexagon in generation $n+2$.
Each pair of $E$-vertices in generation $n$ gives rise to a pair of $F$-vertices across the flat octagon in generation $n+3$.
These are all the $E$-vertices and $F$-vertices, so we have
$$\eqalignno{
E(z) &= {12 z^5 \over 1-z^5} &(2)\cr
F(z) &= {12 z^3 \over 1-z^5}. &(3)\cr
}$$

Let $C(z)$ and $D(z)$ denote the generating functions for $C$-vertices and $D$-vertices, respectively.
Each $A$-vertex in generation $n$ gives rise across the sharp hexagon to a $C$-vertex in generation $n+3$.
Each $B$-vertex in generation $n$ gives rise across the sharp octagon to a $D$-vertex in generation $n+4$.
These are all the $C$-vertices and $D$-vertices, so we have
$$\eqalignno{
C(z) &= z^3 \, A(z) &(4)\cr
D(z) &= z^4 \, B(z). &(5)\cr
}$$

Let $A(z)$ and $B(z)$ denote the generating functions for $A$-vertices and $B$-vertices, respectively.
The one child of each $G$-vertex, $E$-vertex and $C$-vertex is an $A$-vertex.
The one child of each $F$-vertex and $D$-vertex is a $B$-vertex.
One of the two children of each $B$-vertex is an $A$-vertex
(unless it is a $D$-vertex or an $F$-vertex).
One of the two children of each $B$-vertex, and each of the two children of an $A$-vertex,
is a $B$-vertex (unless it is a $C$-vertex, a $D$-vertex or an $E$-vertex).
These observations allow us to write equations for $A(z)$ and $B(z)$
(where the negative terms correspond to the ``unless'' clauses in the observations):
$$\eqalign{
A(z) &= z\,B(z) + z\,C(z) - D(z) + z\,E(z) - F(z) + z\,,G(z) \cr
B(z)&= 2z\,A(z) + z\,B(z) - 2C(z) + z\,D(z) - D(z) - E(z) + z\,F(z). \cr
}$$
Substituting equations (1--5) and solving the resulting equations for $A(z)$ and $B(z)$ yields
$$\eqalign{
A(z) &=
{6z\, (1 - z^2 + 2z^4 + z^5 - z^7) \over (1-z)\,(1+z+z^2+z^3+z^4)\,(1-z^2-2z^3-z^4+z^6)} \cr
B(z) &=
{12z^2\, (1 + z - z^3) \over (1-z)\,(1+z+z^2+z^3+z^4)\,(1-z^2-2z^3-z^4+z^6)}. \cr
}$$
Summing over the seven kinds of vertices yields
$$V(z) = {6 (1+z+z^2+z^3) \over 1-z^2-2z^3-z^4+z^6},$$
as claimed.
\sk

\heading{3. Conclusion}

We have found the generating function that enumerates angels and devils in Escher's {\it Circle Limit IV}.
We did this by considering seven different types of vertices, together with their seven enumerating functions.
The final generating function is strikingly simple as compared with those that appear intermediately in the derivation,
and this of course raises the question of whether there might be some much simpler derivation of this final result.
\sk

\heading{4. Acknowledgment}

The research reported here was supported
by Grant CCF 0646682 from the National Science Foundation.
\sk

\heading{5. References}

\ref B; C. D. Bennett;
``A Paradoxical View of Escher's Angels and Devils'';
Mathematical Intelligencer; 22:3 (2000) 39--46.

\refbook C1; J. Choi;
Counting Vertices in Isohedral Tilings;
B.~S. Thesis, Department of Mathematics,
Harvey Mudd College, Claremont, CA,
May 2012, 27 pp..

\refinbook C2; H. M. S. Coxeter;
``Angels and Devils'';
in: D.~A. Klarner (Ed.);
The Mathematical Gardner;
Springer, 1980, pp.~197--209.

\refinbook D; D. Dunham;
``The Symmetry of ``Circle Limit IV'' and Related Patterns'';
in: C.~S. Kaplan and R.~Sarhangi;
Proceedings of Bridge 2009: 
Mathematics, Music, Art, Architecture, Culture;
pp.~163--168.

\ref P; A. Paul and N. Pippenger;
 ``A Census of Vertices by Generations in Regular Tessellations of the Plane'';
 Electronic Journal of Combinatorics; 18:1 (April 14, 2011) P87 (13 pp.).
 
\bye